\documentclass[reqno]{amsart}

\usepackage{amsmath,amssymb,amsthm,amsfonts,mathrsfs}
\usepackage{fullpage}
\usepackage{xcolor}
\usepackage{graphicx}
\usepackage{listings}
\usepackage{bbm}

\definecolor{codegreen}{rgb}{0,0.6,0}
\definecolor{codegray}{rgb}{0.5,0.5,0.5}
\definecolor{codepurple}{rgb}{0.58,0,0.82}
\definecolor{backcolour}{rgb}{0.95,0.95,0.92}

\lstdefinestyle{mystyle}{
  backgroundcolor=\color{backcolour},   commentstyle=\color{codegreen},
  keywordstyle=\color{magenta},
  numberstyle=\tiny\color{codegray},
  stringstyle=\color{codepurple},
  basicstyle=\ttfamily\footnotesize,
  breakatwhitespace=false,         
  breaklines=true,                 
  captionpos=b,                    
  keepspaces=true,                 
  numbers=left,                    
  numbersep=5pt,                  
  showspaces=false,                
  showstringspaces=false,
  showtabs=false,                  
  tabsize=2
}

\usepackage{listings}
\usepackage[colorlinks=true,
linkcolor=blue,
anchorcolor=blue,
citecolor=red
]{hyperref}
\allowdisplaybreaks % allow equations break between two pages
\newtheorem{theorem}{Theorem}[section]
\newtheorem{lemma}[theorem]{Lemma}

\newtheorem*{conjecture*}{Conjecture}
\theoremstyle{definition}

\theoremstyle{remark}

\newtheorem*{remark*}{remark}
\usepackage{colonequals}

\author{Runbo Li}
\address{International Curriculum Center, The High School Affiliated to Renmin University of China, Beijing, China}
\email{runbo.li.carey@gmail.com}

\makeatletter
\@namedef{subjclassname@2020}{\textup{}2020 Mathematics Subject Classification}
\makeatother

\title[]{Largest square divisors of shifted primes}

\begin{document}
	
\begin{abstract}
The author shows that there are infinitely many primes $p$ such that for any nonzero integer $a$, $p-a$ is divisible by a square $d^2 > p^{\frac{1}{2}+\frac{1}{700}}$. The exponent $\frac{1}{2}+\frac{1}{700}$ improves Merikoski's $\frac{1}{2}+\frac{1}{2000}$. Many powerful devices in Harman's sieve are used for this improvement.
\end{abstract}

\maketitle

\tableofcontents

\section{Introduction}
The Euler's conjecture, which states that there are infinitely many primes of the form $n^2 +1$, is one of Landau's problems on prime numbers. There are several ways to attack this conjecture. One way is to relax the number of prime factors of $f(n)$, and the best result in this way is due to Iwaniec \cite{IwaniecP2}. Building on the previous work of Richert \cite{RichertP3}, he showed that for any irreducible polynomial $f(n)=an^2 +bn+c$ with $a>0$ and $c \equiv 1 (\bmod\ 2)$, there are infinitely many $x$ such that $f(x)$ has at most $2$ prime factors.

Another possible way is to consider the square divisors of $p-1$. If we can show that there are infinitely many primes such that $p-1$ is divisible by a large square $d^2 \geqslant p^{\theta}$ with $\theta = 1$, then the Euler's conjecture is solved. The first result on this direction is due to Baier and Zhao \cite{BZ1} \cite{BZ2}, where they proved the above statement holds with $d$ prime and $\theta < \frac{4}{9}$ as an application of their large sieve for sparse sets of moduli. They interpret the problem as an equidistribution problem for primes $p \equiv 1 (\bmod\ d^2)$, after which the result follows from their Bombieri--Vinogradov type theorem for sparse sets of moduli [\cite{BZ1}, Theorem 3].

In 2009, Matomäki \cite{Matomaki09} improved the above result to $\theta < \frac{1}{2}$ using Harman's sieve \cite{HarmanBOOK} and Type--II information obtained using the large sieve of Baier and Zhao \cite{BZ2}. Note that the exponent $\theta = \frac{1}{2}$ is the limit of what can be obtained under the Generalized Riemann Hypothesis (GRH). In \cite{Merikoski20}, Merikoski first broke this $\frac{1}{2}$--barrier and successfully got $\theta \leqslant \frac{1}{2} + \frac{1}{2000}$ without the restriction that $d$ is a prime. In the article \cite{Merikoski20}, he mentioned that the ``extra'' exponent $\frac{1}{2000}$ has not been fully optimized, and one should be able to increase this to some value between $\frac{1}{500}$ and $\frac{1}{1000}$. In this paper, we increase this to $\frac{1}{700}$ by a careful decomposition on Harman's sieve.

\begin{theorem}\label{t1}
Let $a \neq 0$ be an integer. There are infinitely many primes $p$ such that $d^2 \mid (p - a)$ for some integer $d$ with
$$
d^2 \geqslant p^{\frac{1}{2}+\frac{1}{700}}.
$$
\end{theorem}

Throughout this paper, we always suppose that $\varepsilon$ is a sufficiently small positive constant and $X$ is sufficiently large. The letter $p$, with or without subscript, is reserved for prime numbers. Let $\varpi = \frac{1}{1400}$, $D = X^{\frac{1}{2}+2\varpi}$, $K = \lceil \frac{1}{\varepsilon} \rceil$ and $P = D^{\frac{1}{K}}$. Let $\sigma$ be a number satifsies the condition $19\sigma + 90\varpi + 71\varepsilon <1$. Define
$$
I_j = \left(2^{j-1} P^{\frac{1}{2}}, 2^{j} P^{\frac{1}{2}}\right] \ \text{for} \ j = 1,2,\ldots,K.
$$
We set
$$
\mathcal{D} = \{p_1^2 p_2^2 \cdots p_K^2 : p_j \in I_j \ \text{for} \ j = 1,2,\ldots,K\}
$$
so that $d^2 \in D$ is of size $\asymp D$ and is a square of a squarefree integer.

Fix an integer $a \neq 0$ and a $C^{\infty}$--smooth function $0 \leqslant \psi \leqslant 1$, supported on the interval $[1,2]$ and satisfying $\psi(x) = 1$ for $1+\eta \leqslant x \leqslant 2-\eta$ for some sufficiently small positive $\eta$. For $d^2 \in D$ and $z < X$, denote
$$
S\left(\mathcal{A}^d, z\right) = \sum_{\substack{n \equiv a (\bmod d^2) \\ \left(n, P(z)\right) = 1}} \psi\left(\frac{n}{X}\right) \quad \text{and} \quad S\left(\mathcal{B}^d, z\right) = \frac{1}{\varphi(d^2)} \sum_{\substack{(n,d^2)=1 \\ \left(n, P(z)\right) = 1}} \psi\left(\frac{n}{X}\right).
$$
Then Theorem~\ref{t1} holds if there exists $\varepsilon, \eta, c >0$ such that for all but $O\left(D^{\frac{1}{2}} X^{-\eta}\right)$ of the moduli $d^2 \in \mathcal{D}$, we have
\begin{equation}
S\left(\mathcal{A}^d, 2X^{\frac{1}{2}}\right) \geqslant c S\left(\mathcal{B}^d, 2X^{\frac{1}{2}}\right).
\end{equation}

\section{Asymptotic formulas}

\begin{lemma}\label{l21} ([\cite{Merikoski20}, Proposition 7]).
Let $U \leqslant X^{\frac{1}{2}+2\varpi +\varepsilon}$ and let $a_u$ be divisor--bounded. Then for all but $O\left(D^{\frac{1}{2}} X^{-\eta}\right)$ of $d^2 \in \mathcal{D}$, we have
$$
\sum_{u \sim U} a_u S\left(\mathcal{A}_{u}^d, X^{\sigma - 2\varpi}\right) = (1+o(1)) \sum_{u \sim U} a_u S\left(\mathcal{B}_{u}^d, X^{\sigma - 2\varpi}\right).
$$
\end{lemma}

\begin{lemma}\label{l22} ([\cite{Merikoski20}, Proposition 6]).
Let $U \leqslant X^{\frac{1}{2}-\sigma}$, $V \leqslant X^{\frac{1}{8}+\frac{\sigma}{2}-\frac{5\varpi}{2}-\eta}$ and let $a_u$, $b_v$ be divisor--bounded. Then for all but $O\left(D^{\frac{1}{2}} X^{-\eta}\right)$ of $d^2 \in \mathcal{D}$, we have
$$
\sum_{\substack{u \sim U \\ v \sim V}} a_u b_v S\left(\mathcal{A}_{u v}^d, X^{\sigma - 2\varpi}\right) = (1+o(1)) \sum_{\substack{u \sim U \\ v \sim V}} a_u b_v S\left(\mathcal{B}_{u v}^d, X^{\sigma - 2\varpi}\right).
$$
\end{lemma}

\begin{lemma}\label{l23} ([\cite{Merikoski20}, Proposition 4]).
Let $U V = X$, $X^{\frac{1}{2}-\sigma} \leqslant U \leqslant X^{\frac{1}{2}-2\varpi -\varepsilon}$ and let $a_u$, $b_v$ be divisor--bounded. Then we have
$$
\sum_{d^2 \in \mathcal{D}} \left| \sum_{\substack{u v \equiv a (\bmod d^2) \\ u \sim U, v \sim V}} a_u b_v \psi\left(\frac{uv}{X}\right) - \frac{1}{\varphi(d^2)} \sum_{\substack{(u v,d^2)=1 \\ u \sim U, v \sim V}} a_u b_v \psi\left(\frac{uv}{X}\right) \right| \ll D^{-\frac{1}{2}} X^{1-\eta}.
$$
\end{lemma}

\section{The final decomposition}
Before decomposing, we define asymptotic regions $I$ and $II$ as
\begin{align}
\nonumber I (m, n) :=&\ \left\{ m+n \leqslant \frac{1}{2} +2\varpi, \text{ or } m \leqslant \frac{1}{2} - \sigma \text{ and } n < \frac{1}{8} + \frac{\sigma}{2} -\frac{5\varpi}{2} \right\}, \\
\nonumber II (m, n) :=&\ \left\{ \frac{1}{2} - \sigma \leqslant m \leqslant \frac{1}{2} - 2\varpi \text{ or } \frac{1}{2} - \sigma \leqslant n \leqslant \frac{1}{2} - 2\varpi \right. \\
\nonumber & \left. \qquad \text{or } \frac{1}{2} - \sigma \leqslant m+n \leqslant \frac{1}{2} - 2\varpi \text{ or } \frac{1}{2} + 2\varpi \leqslant m+n \leqslant \frac{1}{2} + \sigma \right\}.
\end{align}
Let $\omega(u)$ denotes the Buchstab function determined by the following differential--difference equation
\begin{align*}
\begin{cases}
\omega(u)=\frac{1}{u}, & \quad 1 \leqslant u \leqslant 2, \\
(u \omega(u))^{\prime}= \omega(u-1), & \quad u \geqslant 2 .
\end{cases}
\end{align*}
Moreover, we have the upper and lower bounds for $\omega(u)$:
\begin{align*}
\omega(u) \geqslant \omega_{0}(u) =
\begin{cases}
\frac{1}{u}, & \quad 1 \leqslant u < 2, \\
\frac{1+\log(u-1)}{u}, & \quad 2 \leqslant u < 3, \\
\frac{1+\log(u-1)}{u} + \frac{1}{u} \int_{2}^{u-1}\frac{\log(t-1)}{t} d t \geqslant 0.5607, & \quad 3 \leqslant u < 4, \\
0.5612, & \quad u \geqslant 4, \\
\end{cases}
\end{align*}
\begin{align*}
\omega(u) \leqslant \omega_{1}(u) =
\begin{cases}
\frac{1}{u}, & \quad 1 \leqslant u < 2, \\
\frac{1+\log(u-1)}{u}, & \quad 2 \leqslant u < 3, \\
\frac{1+\log(u-1)}{u} + \frac{1}{u} \int_{2}^{u-1}\frac{\log(t-1)}{t} d t \leqslant 0.5644, & \quad 3 \leqslant u < 4, \\
0.5617, & \quad u \geqslant 4. \\
\end{cases}
\end{align*}
We shall use $\omega_0(u)$ and $\omega_1(u)$ to give numerical bounds for some sieve functions discussed below. We shall also use the simple upper bound $\omega(u) \leqslant \max(\frac{1}{u}, 0.5672)$ (see Lemma 8(iii) of \cite{JiaPSV}) to estimate high--dimensional integrals. Fix $\sigma = \frac{1}{20.31}$ and let $p_{j}=X^{\alpha_{j}}$. By Buchstab's identity, we have
\begin{align}
\nonumber S\left(\mathcal{A}^d, 2X^{\frac{1}{2}}\right) =&\ S\left(\mathcal{A}^d, X^{\sigma -2\varpi}\right) - \sum_{\sigma -2\varpi \leqslant \alpha_1 < \frac{1}{2}} S\left(\mathcal{A}_{p_1}^d, X^{\sigma -2\varpi}\right) + \sum_{\substack{\sigma -2\varpi \leqslant \alpha_1 < \frac{1}{2} \\ \sigma -2\varpi \leqslant \alpha_2 < \min\left(\alpha_1, \frac{1-\alpha_1}{2}\right) }} S\left(\mathcal{A}_{p_1 p_2}^d, p_2\right) \\
\nonumber =&\ S\left(\mathcal{A}^d, X^{\sigma -2\varpi}\right) - \sum_{\sigma -2\varpi \leqslant \alpha_1 < \frac{1}{2}} S\left(\mathcal{A}_{p_1}^d, X^{\sigma -2\varpi}\right) + \sum_{(\alpha_1, \alpha_2) \in II} S\left(\mathcal{A}_{p_1 p_2}^d, p_2\right) \\
\nonumber & + \sum_{(\alpha_1, \alpha_2) \in A} S\left(\mathcal{A}_{p_1 p_2}^d, p_2\right) + \sum_{(\alpha_1, \alpha_2) \in B} S\left(\mathcal{A}_{p_1 p_2}^d, p_2\right) + \sum_{(\alpha_1, \alpha_2) \in C} S\left(\mathcal{A}_{p_1 p_2}^d, p_2\right) \\
=&\ S_{1} - S_{2} + S_{II} + S_{A} + S_{B} + S_{C},
\end{align}
where
\begin{align}
\nonumber A(\alpha_1, \alpha_2) =&\ \left\{ \sigma -2\varpi \leqslant \alpha_1 < \frac{1}{2},\ \sigma -2\varpi \leqslant \alpha_2 < \min\left(\alpha_1, \frac{1-\alpha_1}{2}\right),\ (\alpha_1, \alpha_2) \notin II, \right. \\
\nonumber & \left. \qquad (\alpha_1, \alpha_2, \alpha_2) \text{ can be partitioned into } (m, n) \in I \right\}, \\
\nonumber B(\alpha_1, \alpha_2) =&\ \left\{ \sigma -2\varpi \leqslant \alpha_1 < \frac{1}{2},\ \sigma -2\varpi \leqslant \alpha_2 < \min\left(\alpha_1, \frac{1-\alpha_1}{2}\right),\ (\alpha_1, \alpha_2) \notin II, \right. \\
\nonumber & \left. \qquad (\alpha_1, \alpha_2, \alpha_2) \text{ cannot be partitioned into } (m, n) \in I,\ (\alpha_1, \alpha_2) \in I,\ (1-\alpha_1-\alpha_2, \alpha_2) \in I \right\}, \\
\nonumber C(\alpha_1, \alpha_2) =&\ \left\{ \sigma -2\varpi \leqslant \alpha_1 < \frac{1}{2},\ \sigma -2\varpi \leqslant \alpha_2 < \min\left(\alpha_1, \frac{1-\alpha_1}{2}\right),\ (\alpha_1, \alpha_2) \notin II \cup A \cup B \right\}.
\end{align}
We have asymptotic formulas for $S_{1}$ and $S_{2}$ by Lemma~\ref{l21}. We can also give an asymptotic formula for $S_{II}$ by Lemma~\ref{l23}. For $S_{A}$, we can apply Buchstab's identity to get
\begin{align}
\nonumber S_{A} =&\ \sum_{(\alpha_1, \alpha_2) \in A} S\left(\mathcal{A}_{p_1 p_2}^d, p_2\right) \\
\nonumber =&\ \sum_{(\alpha_1, \alpha_2) \in A} S\left(\mathcal{A}_{p_1 p_2}^d, X^{\sigma -2\varpi}\right) - \sum_{\substack{(\alpha_1, \alpha_2) \in A \\ \sigma -2\varpi \leqslant \alpha_3 < \min\left(\alpha_2, \frac{1-\alpha_1-\alpha_2}{2}\right) \\ (\alpha_1, \alpha_2, \alpha_3) \text{ can be partitioned into } (m, n) \in II }} S\left(\mathcal{A}_{p_1 p_2 p_3}^d, p_3\right) \\
\nonumber & - \sum_{\substack{(\alpha_1, \alpha_2) \in A \\ \sigma -2\varpi \leqslant \alpha_3 < \min\left(\alpha_2, \frac{1-\alpha_1-\alpha_2}{2}\right) \\ (\alpha_1, \alpha_2, \alpha_3) \text{ cannot be partitioned into } (m, n) \in II }} S\left(\mathcal{A}_{p_1 p_2 p_3}^d, X^{\sigma -2\varpi}\right) \\
\nonumber & + \sum_{\substack{(\alpha_1, \alpha_2) \in A \\ \sigma -2\varpi \leqslant \alpha_3 < \min\left(\alpha_2, \frac{1-\alpha_1-\alpha_2}{2}\right) \\ (\alpha_1, \alpha_2, \alpha_3) \text{ cannot be partitioned into } (m, n) \in II \\ \sigma -2\varpi \leqslant \alpha_4 < \min\left(\alpha_3, \frac{1-\alpha_1-\alpha_2-\alpha_3}{2}\right) \\ (\alpha_1, \alpha_2, \alpha_3, \alpha_4) \text{ can be partitioned into } (m, n) \in II }} S\left(\mathcal{A}_{p_1 p_2 p_3 p_4}^d, p_4\right) \\
\nonumber & + \sum_{\substack{(\alpha_1, \alpha_2) \in A \\ \sigma -2\varpi \leqslant \alpha_3 < \min\left(\alpha_2, \frac{1-\alpha_1-\alpha_2}{2}\right) \\ (\alpha_1, \alpha_2, \alpha_3) \text{ cannot be partitioned into } (m, n) \in II \\ \sigma -2\varpi \leqslant \alpha_4 < \min\left(\alpha_3, \frac{1-\alpha_1-\alpha_2-\alpha_3}{2}\right) \\ (\alpha_1, \alpha_2, \alpha_3, \alpha_4) \text{ cannot be partitioned into } (m, n) \in II }} S\left(\mathcal{A}_{p_1 p_2 p_3 p_4}^d, p_4\right) \\
=&\ S_{A1} - S_{A2} - S_{A3} + S_{A4} + S_{A5}.
\end{align}
By Lemmas~\ref{l21}--\ref{l22} we have asymptotic formulas for $S_{A1}$ and $S_{A3}$. For $S_{A2}$ and $S_{A4}$ we can also give asymptotic formulas by Lemma~\ref{l23}. We discard part of $S_{A5}$ if we can neither give an asymptotic formula nor decompose it further. For the remaining part, we can perform Buchstab's identity twice more to reach a six--dimensional sum if we can group $(\alpha_1, \alpha_2, \alpha_3, \alpha_4, \alpha_4)$ into $(m, n) \in I$, and we can use reversed Buchstab's identity to make some almost--primes visible. Working as in \cite{LRBapmod1} and \cite{LRB052}, the total loss from $S_{A}$ can be bounded by
\begin{align}
\nonumber & \left( \int_{(t_1, t_2, t_3, t_4) \in S_{A51}} \frac{\omega_1 \left(\frac{1 - t_1 - t_2 - t_3 - t_4}{t_4}\right)}{t_1 t_2 t_3 t_4^2} d t_4 d t_3 d t_2 d t_1 \right) \\
\nonumber -& \left( \int_{(t_1, t_2, t_3, t_4, t_5) \in S_{A52}} \frac{\omega_0 \left(\frac{1 - t_1 - t_2 - t_3 - t_4 - t_5}{t_5}\right)}{t_1 t_2 t_3 t_4 t_5^2} d t_5 d t_4 d t_3 d t_2 d t_1 \right) \\
\nonumber +& \left( \int_{(t_1, t_2, t_3, t_4, t_5, t_6) \in S_{A53}} \frac{\omega_1 \left(\frac{1 - t_1 - t_2 - t_3 - t_4 - t_5 - t_6}{t_6}\right)}{t_1 t_2 t_3 t_4 t_5 t_6^2} d t_6 d t_5 d t_4 d t_3 d t_2 d t_1 \right) \\
\nonumber +& \left( \int_{(t_1, t_2, t_3, t_4, t_5, t_6, t_7, t_8) \in S_{A54}} \frac{\max \left(\frac{t_8}{1 - t_1 - t_2 - t_3 - t_4 - t_5 - t_6 - t_7 - t_8}, 0.5672\right)}{t_1 t_2 t_3 t_4 t_5 t_6 t_7 t_8^2} d t_8 d t_7 d t_6 d t_5 d t_4 d t_3 d t_2 d t_1 \right) \\
<&\ 0.002515,
\end{align}
where
\begin{align}
\nonumber S_{A51}(t_1, t_2, t_3, t_4) :=&\ \left\{ (t_1, t_2) \in S_{A}, \right. \\
\nonumber & \quad \sigma -2\varpi \leqslant t_3 < \min\left(t_2, \frac{1}{2}(1-t_1-t_2)\right), \\ 
\nonumber & \quad (t_1, t_2, t_3) \text{ cannot be partitioned into } (m, n) \in II, \\
\nonumber & \quad \sigma -2\varpi \leqslant t_4 < \min\left(t_3, \frac{1}{2}(1-t_1-t_2-t_3) \right), \\
\nonumber & \quad (t_1, t_2, t_3, t_4) \text{ cannot be partitioned into } (m, n) \in II, \\
\nonumber & \quad (t_1, t_2, t_3, t_4, t_4) \text{ cannot be partitioned into } (m, n) \in I, \\
\nonumber & \left. \quad \sigma -2\varpi \leqslant t_1 < \frac{1}{2},\ \sigma -2\varpi \leqslant t_2 < \min\left(t_1, \frac{1}{2}(1-t_1) \right) \right\}, \\
\nonumber S_{A52}(t_1, t_2, t_3, t_4, t_5) :=&\ \left\{ (t_1, t_2) \in S_{A}, \right. \\
\nonumber & \quad \sigma -2\varpi \leqslant t_3 < \min\left(t_2, \frac{1}{2}(1-t_1-t_2)\right), \\ 
\nonumber & \quad (t_1, t_2, t_3) \text{ cannot be partitioned into } (m, n) \in II, \\
\nonumber & \quad \sigma -2\varpi \leqslant t_4 < \min\left(t_3, \frac{1}{2}(1-t_1-t_2-t_3) \right), \\
\nonumber & \quad (t_1, t_2, t_3, t_4) \text{ cannot be partitioned into } (m, n) \in II, \\
\nonumber & \quad (t_1, t_2, t_3, t_4, t_4) \text{ cannot be partitioned into } (m, n) \in I, \\
\nonumber & \quad t_4 < t_5 < \frac{1}{2}(1-t_1-t_2-t_3-t_4), \\
\nonumber & \quad (t_1, t_2, t_3, t_4, t_5) \text{ can be partitioned into } (m, n) \in II, \\
\nonumber & \left. \quad \sigma -2\varpi \leqslant t_1 < \frac{1}{2},\ \sigma -2\varpi \leqslant t_2 < \min\left(t_1, \frac{1}{2}(1-t_1) \right) \right\}, \\
\nonumber S_{A53}(t_1, t_2, t_3, t_4, t_5, t_6) :=&\ \left\{ (t_1, t_2) \in S_{A}, \right. \\
\nonumber & \quad \sigma -2\varpi \leqslant t_3 < \min\left(t_2, \frac{1}{2}(1-t_1-t_2)\right), \\ 
\nonumber & \quad (t_1, t_2, t_3) \text{ cannot be partitioned into } (m, n) \in II, \\
\nonumber & \quad \sigma -2\varpi \leqslant t_4 < \min\left(t_3, \frac{1}{2}(1-t_1-t_2-t_3) \right), \\
\nonumber & \quad (t_1, t_2, t_3, t_4) \text{ cannot be partitioned into } (m, n) \in II, \\
\nonumber & \quad (t_1, t_2, t_3, t_4, t_4) \text{ can be partitioned into } (m, n) \in I, \\
\nonumber & \quad \sigma -2\varpi \leqslant t_5 < \min\left(t_4, \frac{1}{2}(1-t_1-t_2-t_3-t_4) \right), \\
\nonumber & \quad (t_1, t_2, t_3, t_4, t_5) \text{ cannot be partitioned into } (m, n) \in II, \\
\nonumber & \quad \sigma -2\varpi \leqslant t_6 < \min\left(t_5, \frac{1}{2}(1-t_1-t_2-t_3-t_4-t_5) \right),\\
\nonumber & \quad (t_1, t_2, t_3, t_4, t_5, t_6) \text{ cannot be partitioned into } (m, n) \in II, \\
\nonumber & \quad (t_1, t_2, t_3, t_4, t_5, t_6, t_6) \text{ cannot be partitioned into } (m, n) \in I, \\
\nonumber & \left. \quad \sigma -2\varpi \leqslant t_1 < \frac{1}{2},\ \sigma -2\varpi \leqslant t_2 < \min\left(t_1, \frac{1}{2}(1-t_1) \right) \right\}, \\
\nonumber S_{A54}(t_1, t_2, t_3, t_4, t_5, t_6, t_7, t_8) :=&\ \left\{ (t_1, t_2) \in S_{A}, \right. \\
\nonumber & \quad \sigma -2\varpi \leqslant t_3 < \min\left(t_2, \frac{1}{2}(1-t_1-t_2)\right), \\ 
\nonumber & \quad (t_1, t_2, t_3) \text{ cannot be partitioned into } (m, n) \in II, \\
\nonumber & \quad \sigma -2\varpi \leqslant t_4 < \min\left(t_3, \frac{1}{2}(1-t_1-t_2-t_3) \right), \\
\nonumber & \quad (t_1, t_2, t_3, t_4) \text{ cannot be partitioned into } (m, n) \in II, \\
\nonumber & \quad (t_1, t_2, t_3, t_4, t_4) \text{ can be partitioned into } (m, n) \in I, \\
\nonumber & \quad \sigma -2\varpi \leqslant t_5 < \min\left(t_4, \frac{1}{2}(1-t_1-t_2-t_3-t_4) \right), \\
\nonumber & \quad (t_1, t_2, t_3, t_4, t_5) \text{ cannot be partitioned into } (m, n) \in II, \\
\nonumber & \quad \sigma -2\varpi \leqslant t_6 < \min\left(t_5, \frac{1}{2}(1-t_1-t_2-t_3-t_4-t_5) \right),\\
\nonumber & \quad (t_1, t_2, t_3, t_4, t_5, t_6) \text{ cannot be partitioned into } (m, n) \in II, \\
\nonumber & \quad (t_1, t_2, t_3, t_4, t_5, t_6, t_6) \text{ can be partitioned into } (m, n) \in I, \\
\nonumber & \quad \sigma -2\varpi \leqslant t_7 < \min\left(t_5, \frac{1}{2}(1-t_1-t_2-t_3-t_4-t_5-t_6) \right),\\
\nonumber & \quad (t_1, t_2, t_3, t_4, t_5, t_6, t_7) \text{ cannot be partitioned into } (m, n) \in II, \\
\nonumber & \quad \sigma -2\varpi \leqslant t_8 < \min\left(t_5, \frac{1}{2}(1-t_1-t_2-t_3-t_4-t_5-t_6-t_7) \right),\\
\nonumber & \quad (t_1, t_2, t_3, t_4, t_5, t_6, t_7, t_8) \text{ cannot be partitioned into } (m, n) \in II, \\
\nonumber & \left. \quad \sigma -2\varpi \leqslant t_1 < \frac{1}{2},\ \sigma -2\varpi \leqslant t_2 < \min\left(t_1, \frac{1}{2}(1-t_1) \right) \right\}.
\end{align}

For $S_{B}$ we use Buchstab's identity to get
\begin{equation}
S_B = \sum_{(\alpha_1, \alpha_2) \in B} S\left(\mathcal{A}_{p_1 p_2}^d, X^{\sigma -2\varpi}\right) - \sum_{\substack{(\alpha_1, \alpha_2) \in B \\ \sigma -2\varpi \leqslant \alpha_3 < \min\left(\alpha_2, \frac{1-\alpha_1-\alpha_2}{2}\right) }} S\left(\mathcal{A}_{p_1 p_2 p_3}^d, p_3\right).
\end{equation}
Here we cannot decompose it directly using Buchstab's identity once more, since we cannot give an asymptotic formula for part of the negative sum
$$
\sum_{\substack{(\alpha_1, \alpha_2) \in B \\ \sigma -2\varpi \leqslant \alpha_3 < \min\left(\alpha_2, \frac{1-\alpha_1-\alpha_2}{2}\right) \\ (\alpha_1, \alpha_2, \alpha_3) \text{ cannot be partitioned into } (m, n) \in II }} S\left(\mathcal{A}_{p_1 p_2 p_3}^d, X^{\sigma -2\varpi}\right).
$$
However, we can use a role--reversal to transfer the last sum in (5) into a form that have an asymptotic formula. Note that in \cite{Merikoski20} role--reversals were not used. We refer the readers to \cite{LiRunbo1215} and \cite{LRB052} for more applications of role--reversals. By a standard process, we have (where $\beta \sim X^{1-\alpha_1 -\alpha_2 -\alpha_3}$ and $\left(\beta, P(p_{3})\right)=1$)
\begin{align}
\nonumber S_{B} =&\ \sum_{(\alpha_1, \alpha_2) \in B} S\left(\mathcal{A}_{p_1 p_2}^d, p_2\right) \\
\nonumber =&\ \sum_{(\alpha_1, \alpha_2) \in B} S\left(\mathcal{A}_{p_1 p_2}^d, X^{\sigma -2\varpi}\right) - \sum_{\substack{(\alpha_1, \alpha_2) \in B \\ \sigma -2\varpi \leqslant \alpha_3 < \min\left(\alpha_2, \frac{1-\alpha_1-\alpha_2}{2}\right) \\ (\alpha_1, \alpha_2, \alpha_3) \text{ can be partitioned into } (m, n) \in II }} S\left(\mathcal{A}_{p_1 p_2 p_3}^d, p_3\right) \\
\nonumber & - \sum_{\substack{(\alpha_1, \alpha_2) \in B \\ \sigma -2\varpi \leqslant \alpha_3 < \min\left(\alpha_2, \frac{1-\alpha_1-\alpha_2}{2}\right) \\ (\alpha_1, \alpha_2, \alpha_3) \text{ cannot be partitioned into } (m, n) \in II }} S\left(\mathcal{A}_{\beta p_2 p_3}^d, X^{\sigma -2\varpi}\right) \\
\nonumber & + \sum_{\substack{(\alpha_1, \alpha_2) \in B \\ \sigma -2\varpi \leqslant \alpha_3 < \min\left(\alpha_2, \frac{1-\alpha_1-\alpha_2}{2}\right) \\ (\alpha_1, \alpha_2, \alpha_3) \text{ cannot be partitioned into } (m, n) \in II \\ \sigma -2\varpi \leqslant \alpha_4 < \frac{1}{2}\alpha_1 \\ (1-\alpha_1-\alpha_2-\alpha_3, \alpha_2, \alpha_3, \alpha_4) \text{ can be partitioned into } (m, n) \in II }} S\left(\mathcal{A}_{\beta p_2 p_3 p_4}^d, p_4\right) \\
\nonumber & + \sum_{\substack{(\alpha_1, \alpha_2) \in B \\ \sigma -2\varpi \leqslant \alpha_3 < \min\left(\alpha_2, \frac{1-\alpha_1-\alpha_2}{2}\right) \\ (\alpha_1, \alpha_2, \alpha_3) \text{ cannot be partitioned into } (m, n) \in II \\ \sigma -2\varpi \leqslant \alpha_4 < \frac{1}{2}\alpha_1 \\ (1-\alpha_1-\alpha_2-\alpha_3, \alpha_2, \alpha_3, \alpha_4) \text{ cannot be partitioned into } (m, n) \in II }} S\left(\mathcal{A}_{\beta p_2 p_3 p_4}^d, p_4\right) \\
=&\ S_{B1} - S_{B2} - S_{B3} + S_{B4} + S_{B5}.
\end{align}
We can give asymptotic formulas for $S_{B1}$--$S_{B4}$ by Lemmas~\ref{l21}--\ref{l23}. We can also decompose part of $S_{B5}$ if the variables (with $\alpha_1$ replaced by $1-\alpha_1-\alpha_2-\alpha_3$) satisfy the same conditions as in the decomposable part of $S_{B5}$. Again, the total loss from $S_{B}$ can be bounded by
\begin{align}
\nonumber & \left( \int_{(t_1, t_2, t_3, t_4) \in S_{B51}} \frac{\omega_1 \left(\frac{t_1 - t_4}{t_4}\right) \omega_1 \left(\frac{1 - t_1 - t_2 - t_3}{t_3}\right)}{t_2 t_3^2 t_4^2} d t_4 d t_3 d t_2 d t_1 \right) \\
\nonumber -& \left( \int_{(t_1, t_2, t_3, t_4, t_5) \in S_{B52}} \frac{\omega_0 \left(\frac{t_1 - t_4 - t_5}{t_5}\right) \omega_0 \left(\frac{1 - t_1 - t_2 - t_3}{t_3}\right)}{t_2 t_3^2 t_4 t_5^2} d t_5 d t_4 d t_3 d t_2 d t_1 \right) \\
\nonumber +& \left( \int_{(t_1, t_2, t_3, t_4, t_5, t_6) \in S_{B53}} \frac{\omega_1 \left(\frac{t_1 - t_4 - t_5 - t_6}{t_6}\right) \omega_1 \left(\frac{1 - t_1 - t_2 - t_3}{t_3}\right)}{t_2 t_3^2 t_4 t_5 t_6^2} d t_6 d t_5 d t_4 d t_3 d t_2 d t_1 \right) \\
<&\ 0.006249,
\end{align}
where
\begin{align}
\nonumber S_{B51}(t_1, t_2, t_3, t_4) :=&\ \left\{ (t_1, t_2) \in S_{B}, \right. \\
\nonumber & \quad \sigma -2\varpi \leqslant t_3 < \min\left(t_2, \frac{1}{2}(1-t_1-t_2)\right), \\ 
\nonumber & \quad (t_1, t_2, t_3) \text{ cannot be partitioned into } (m, n) \in II, \\
\nonumber & \quad \sigma -2\varpi \leqslant t_4 < \frac{1}{2}t_1, \\
\nonumber & \quad (1-t_1-t_2-t_3, t_2, t_3, t_4) \text{ cannot be partitioned into } (m, n) \in II, \\
\nonumber & \quad (1-t_1-t_2-t_3, t_2, t_3, t_4, t_4) \text{ cannot be partitioned into } (m, n) \in I, \\
\nonumber & \left. \quad \sigma -2\varpi \leqslant t_1 < \frac{1}{2},\ \sigma -2\varpi \leqslant t_2 < \min\left(t_1, \frac{1}{2}(1-t_1) \right) \right\}, \\
\nonumber S_{B52}(t_1, t_2, t_3, t_4, t_5) :=&\ \left\{ (t_1, t_2) \in S_{B}, \right. \\
\nonumber & \quad \sigma -2\varpi \leqslant t_3 < \min\left(t_2, \frac{1}{2}(1-t_1-t_2)\right), \\ 
\nonumber & \quad (t_1, t_2, t_3) \text{ cannot be partitioned into } (m, n) \in II, \\
\nonumber & \quad \sigma -2\varpi \leqslant t_4 < \frac{1}{2}t_1, \\
\nonumber & \quad (1-t_1-t_2-t_3, t_2, t_3, t_4) \text{ cannot be partitioned into } (m, n) \in II, \\
\nonumber & \quad (1-t_1-t_2-t_3, t_2, t_3, t_4, t_4) \text{ cannot be partitioned into } (m, n) \in I, \\
\nonumber & \quad t_4 < t_5 < \frac{1}{2}(t_1-t_4), \\
\nonumber & \quad (1-t_1-t_2-t_3, t_2, t_3, t_4, t_5) \text{ can be partitioned into } (m, n) \in II, \\
\nonumber & \left. \quad \sigma -2\varpi \leqslant t_1 < \frac{1}{2},\ \sigma -2\varpi \leqslant t_2 < \min\left(t_1, \frac{1}{2}(1-t_1) \right) \right\}, \\
\nonumber S_{B53}(t_1, t_2, t_3, t_4, t_5, t_6) :=&\ \left\{ (t_1, t_2) \in S_{B}, \right. \\
\nonumber & \quad \sigma -2\varpi \leqslant t_3 < \min\left(t_2, \frac{1}{2}(1-t_1-t_2)\right), \\ 
\nonumber & \quad (t_1, t_2, t_3) \text{ cannot be partitioned into } (m, n) \in II, \\
\nonumber & \quad \sigma -2\varpi \leqslant t_4 < \frac{1}{2}t_1, \\
\nonumber & \quad (1-t_1-t_2-t_3, t_2, t_3, t_4) \text{ cannot be partitioned into } (m, n) \in II, \\
\nonumber & \quad (1-t_1-t_2-t_3, t_2, t_3, t_4, t_4) \text{ can be partitioned into } (m, n) \in I, \\
\nonumber & \quad \sigma -2\varpi \leqslant t_5 < \min\left(t_4, \frac{1}{2}(t_1-t_4) \right), \\
\nonumber & \quad (1-t_1-t_2-t_3, t_2, t_3, t_4, t_5) \text{ cannot be partitioned into } (m, n) \in II, \\
\nonumber & \quad \sigma -2\varpi \leqslant t_6 < \min\left(t_5, \frac{1}{2}(t_1-t_4-t_5) \right),\\
\nonumber & \quad (1-t_1-t_2-t_3, t_2, t_3, t_4, t_5, t_6) \text{ cannot be partitioned into } (m, n) \in II, \\
\nonumber & \left. \quad \sigma -2\varpi \leqslant t_1 < \frac{1}{2},\ \sigma -2\varpi \leqslant t_2 < \min\left(t_1, \frac{1}{2}(1-t_1) \right) \right\}.
\end{align}

For $S_{C}$ we cannot use Buchstab's identity in a straightforward manner, but we can use Buchstab's identity in reverse to make some almost--primes visible. The details of using Buchstab's identity in reverse are similar to those in \cite{LiRunbo1215} and \cite{LRB052}. By using Buchstab's identity in reverse twice, we have
\begin{align}
\nonumber S_{C} =&\ \sum_{(\alpha_1, \alpha_2) \in C} S\left(\mathcal{A}_{p_1 p_2}^d, p_2\right) \\
\nonumber =&\ \sum_{(\alpha_1, \alpha_2) \in C} S\left(\mathcal{A}_{p_1 p_2}^d, 2\left(\frac{X}{p_1 p_2}\right)^{\frac{1}{2}} \right) \\
\nonumber &+ \sum_{\substack{(\alpha_1, \alpha_2) \in C \\ \alpha_2 < \alpha_3 < \frac{1-\alpha_1-\alpha_2}{2} \\ (\alpha_1, \alpha_2, \alpha_3) \text{ can be partitioned into } (m, n) \in II }} S\left(\mathcal{A}_{p_1 p_2 p_3}^d, p_3\right)\\
\nonumber &+ \sum_{\substack{(\alpha_1, \alpha_2) \in C \\ \alpha_2 < \alpha_3 < \frac{1-\alpha_1-\alpha_2}{2} \\ (\alpha_1, \alpha_2, \alpha_3) \text{ cannot be partitioned into } (m, n) \in II }} S\left(\mathcal{A}_{p_1 p_2 p_3}^d, 2\left(\frac{X}{p_1 p_2 p_3}\right)^{\frac{1}{2}} \right)\\
\nonumber &+ \sum_{\substack{(\alpha_1, \alpha_2) \in C \\ \alpha_2 < \alpha_3 < \frac{1-\alpha_1-\alpha_2}{2} \\ (\alpha_1, \alpha_2, \alpha_3) \text{ cannot be partitioned into } (m, n) \in II \\ \alpha_3 < \alpha_4 < \frac{1-\alpha_1-\alpha_2-\alpha_3}{2} \\ (\alpha_1, \alpha_2, \alpha_3, \alpha_4) \text{ can be partitioned into } (m, n) \in II }} S\left(\mathcal{A}_{p_1 p_2 p_3 p_4}^d, p_4\right)\\
\nonumber &+ \sum_{\substack{(\alpha_1, \alpha_2) \in C \\ \alpha_2 < \alpha_3 < \frac{1-\alpha_1-\alpha_2}{2} \\ (\alpha_1, \alpha_2, \alpha_3) \text{ cannot be partitioned into } (m, n) \in II \\ \alpha_3 < \alpha_4 < \frac{1-\alpha_1-\alpha_2-\alpha_3}{2} \\ (\alpha_1, \alpha_2, \alpha_3, \alpha_4) \text{ cannot be partitioned into } (m, n) \in II }} S\left(\mathcal{A}_{p_1 p_2 p_3 p_4}^d, p_4\right)\\
=&\ S_{C1} + S_{C2} + S_{C3} + S_{C4} + S_{C5}.
\end{align}
We can give asymptotic formulas for $S_{C2}$ and $S_{C4}$ by Lemma~\ref{l23}, hence we can subtract them from the loss. In this way we obtain a loss from $S_{C}$ of
\begin{align}
\nonumber & \left( \int_{(t_1, t_2) \in S_C} \frac{\omega\left(\frac{1 - t_1 - t_2}{t_2}\right)}{t_1 t_2^2} d t_2 d t_1 \right) \\
\nonumber -& \left( \int_{(t_1, t_2, t_3) \in S_{C2}} \frac{\omega\left(\frac{1 - t_1 - t_2 - t_3}{t_3}\right)}{t_1 t_2 t_3^2} d t_3 d t_2 d t_1 \right) \\
\nonumber -& \left( \int_{(t_1, t_2, t_3, t_4) \in S_{C4}} \frac{\omega\left(\frac{1 - t_1 - t_2 - t_3 - t_4}{t_4}\right)}{t_1 t_2 t_3 t_4^2} d t_4 d t_3 d t_2 d t_1 \right) \\
<&\ 0.990258,
\end{align}
where
\begin{align}
\nonumber S_{C2}(t_1, t_2, t_3) :=&\ \left\{ (t_1, t_2) \in S_{C},\ t_2 < t_3 < \frac{1}{2}(1-t_1-t_2), \right. \\ 
\nonumber & \quad (t_1, t_2, t_3) \text{ can be partitioned into } (m, n) \in II, \\
\nonumber & \left. \quad \sigma -2\varpi \leqslant t_1 < \frac{1}{2},\ \sigma -2\varpi \leqslant t_2 < \min\left(t_1, \frac{1}{2}(1-t_1) \right) \right\}, \\
\nonumber S_{C4}(t_1, t_2, t_3, t_4) :=&\ \left\{ (t_1, t_2) \in S_{C},\ t_2 < t_3 < \frac{1}{2}(1-t_1-t_2), \right. \\ 
\nonumber & \quad (t_1, t_2, t_3) \text{ cannot be partitioned into } (m, n) \in II, \\
\nonumber & \quad t_3 < t_4 < \frac{1}{2}(1-t_1-t_2-t_3), \\
\nonumber & \quad (t_1, t_2, t_3, t_4) \text{ can be partitioned into } (m, n) \in II, \\
\nonumber & \left. \quad \sigma -2\varpi \leqslant t_1 < \frac{1}{2},\ \sigma -2\varpi \leqslant t_2 < \min\left(t_1, \frac{1}{2}(1-t_1) \right) \right\}.
\end{align}

Finally, by (2)--(9), the total loss is less than
$$
0.990258 + 0.002515 + 0.006249 < 0.9991 < 1
$$
and the proof of Theorem~\ref{t1} is completed.

\bibliographystyle{plain}
\bibliography{bib}

\end{document}